\documentstyle[twoside,epic,eepic,12pt]{article}
\title{\bf WHEN SCHREIER TRANSVERSALS GROW WILD}

\author
{Amnon Rosenmann \\
\vspace {-2 mm}
\small
Dept. of Math. \& Computer Science \\
\vspace {-2 mm}
\small
Ben-Gurion University \\
\vspace {-2 mm}
\small
Beer-Sheva, Israel \\
\vspace {-2 mm}
\small
aro@black.bgu.ac.il
}

\date{}

\begin{document}
\maketitle

\newtheorem{claim}{Claim}[section]
\newtheorem{remark}[claim]{Remark}
\newtheorem{example}[claim]{Example}
\newtheorem{definition}[claim]{Definition}
\newtheorem{theorem}[claim]{Theorem}
\newtheorem{lemma}[claim]{Lemma}
\newtheorem{proposition}[claim]{Proposition}
\newtheorem{corollary}[claim]{Corollary}
\newtheorem{conjecture}[claim]{Conjecture}

\newcommand{\mathcal} {\cal}
\newcommand{\CC} {\mbox{$\mathcal C$}}
\newcommand{\DD} {\mbox{$\mathcal D$}}
\newcommand{\EE} {\mbox{$\mathcal E$}}
\newcommand{\CS} {\mbox{$\mathcal S$}}
\newcommand{\TT} {\mbox{$\mathcal T$}}
\newcommand{\dout} {\mbox{$\scriptstyle out$}}
\newcommand{\ds} {\mbox{$\scriptscriptstyle S$}}
\newcommand{\dx} {\mbox{$\scriptscriptstyle X$}}
\newcommand{\dt} {\mbox{$\scriptscriptstyle \mathcal{T}$}}
\newcommand{\dd} {\mbox{$\scriptscriptstyle D$}}
\newcommand{\bdt} {\mbox{$\scriptscriptstyle \overline{\mathcal T}$}}
\newcommand{\dfhB} {\mbox{$\scriptscriptstyle F/H_2$}}
\newcommand{\dfhA} {\mbox{$\scriptscriptstyle F/H_1$}}
\newcommand{\dhB} {\mbox{$\scriptscriptstyle H_2$}}
\newcommand{\dhA} {\mbox{$\scriptscriptstyle H_1$}}

\newcommand{\dfh} {\mbox{$\scriptscriptstyle F/H$}}
\newcommand{\dg} {\mbox{$\scriptscriptstyle G$}}
\newcommand{\dh} {\mbox{$\scriptscriptstyle H$}}
\newcommand{\df} {\mbox{$\scriptscriptstyle F$}}
\newcommand{\dfa} {\mbox{$\scriptscriptstyle F/A$}}
\newcommand{\da} {\mbox{$\scriptscriptstyle A$}}
\newcommand{\dc} {\mbox{$\scriptscriptstyle C$}}

\newcommand{\lrw}{\mbox{$\longrightarrow$}}
\newcommand{\Lrw}{\mbox{$\Longrightarrow$}}
\newcommand{\Llrw}{\mbox{$\Longleftrightarrow$}}


\baselineskip=12pt  
\baselineskip=6pt  
\small
\section{Introduction}
Schreier formula for the rank of a subgroup of finite index
of a finitely generated free group $F$ is generalized to an arbitrary
(even infinitely generated) subgroup $H$ through the Schreier transversals
of $H$ in $F$. The rank formula may also be expressed in terms of the
cogrowth of $H$. \\

We introduce the rank-growth function $rk_{\dh}(i)$ of a subgroup $H$
of a finitely generated free group $F$. $rk_{\dh}(i)$ is defined to be
the rank of the subgroup of $H$ generated by elements of length
less than or equal to $i$ (with respect to the generators of $F$),
and it equals the
rank of the fundamental group of the subgraph of the cosets graph of $H$,
which consists of the paths starting at $1$  that are of length $\leq i$.
When $H$ is supnormal, i.e.
contains a non-trivial normal subgroup of $F$, we show that its rank-growth
is equivalent to the cogrowth of $H$. A special case of this is the known
result that a supnormal subgroup of $F$ is of finite index if and only
if it is finitely generated. In particular, when $H$ is normal then the
growth of the group $G=F/H$ is equivalent to the rank-growth of $H$. \\

A Schreier transversal forms a spanning tree of the cosets
graph of $H$, and thus its topological structure is of a contractible
spanning subcomplex of a simplicial complex. The $d$-dimensional
simplicial complexes that contain contractible spanning subcomplexes
have the homotopy type of a bouquet of $r$ $d$-spheres.
When these complexes are also $n$-regular then $r$ can be computed by
generalizing the rank formula (which applies to Schreier transversals)
to higher dimensions. \\

Let us remark that part of the results here apply in a similar form
also to Schreier transversals and Schreier bases of right ideals in
free group algebras (see \cite{Lew}, \cite{RR}, \cite{Ros2}).

\section{Generalized Schreier Formula}
Let $H$ be a non-trivial subgroup of a free group $F$.
By the Nielsen-Schreier Theorem $H$ is free too
(see, for example, \cite{Lyn}), and explicit free generators for
it can be given.
Suppose that $F$ is freely generated on a set $X$ (not necessarily
finite). The Cayley graph of $F$ (with respect to $X$) has the form of
a tree and is the universal covering of a space $\mathcal Q$ which
is a bouquet of $|X|$ loops ($| . |$ denotes cardinality throughout the paper).
The covering space of $\mathcal Q$ with regard to $H$ is
the {\em cosets graph} $\mathcal G$ of $H$, and it is obtained as the
quotient of the Cayley graph of $F$ under the left action of $H$.
Thus $H$ is the fundamental group of $\mathcal G$.
The set of vertices of $\mathcal G$ is the set of right
cosets of $H$ in $F$. A (double) edge which is labeled with $x \in X$
goes in the direction from the coset $Ht_1$ to the coset $Ht_2$ if
and only if $Ht_2=Ht_1x$, and it is labeled with
$x^{-1}$ in the direction from $Ht_2$ to $Ht_1$. This gives a
connected graph with $|X \cup X^{-1}|$ edges at each vertex.
It is more convenient to
label the vertices of $\mathcal G$ with specific coset representatives
in the following way.  Let $\mathcal T$ be a spanning tree of $\mathcal G$.
The identity element $1$
is chosen to represent the coset $H$ and defined to be the root
of $\mathcal T$, and each other vertex is labeled with
the group element one gets by reading off the edge labels in
a path in $\mathcal T$ that starts at the root and ends at the given vertex.
We also denote by $\mathcal T$ the set of (the labels of) the
vertices $V({\mathcal T})$
of the tree $\mathcal T$, that is the coset representatives of $H$.
This set is a {\em Schreier transversal} for $H$ in $F$, which is
characterized by the property that every initial segment
of an element of $\mathcal T$ is also in $\mathcal T$.
For each $1 \neq w \in H$
there exist $u,v \in \mathcal T$ of maximal lengths such that $u$ is a prefix
of $w$ and $v$ is a prefix of $w^{-1}$. Since $t_1t_2^{-1} \notin H$ for
every $t_1 \neq t_2$ in $\mathcal T$, then $l(u)+l(v)<l(w)$, where $l$ denotes
the length of the (reduced) element in $F$. The Schreier generators
for $H$ with respect to $\mathcal T$ are those $w \in H$ for which
\begin{equation}
l(u)+l(v)= l(w)-1.
\end{equation}
Moreover, if $\phi$ is the coset map associated with $\mathcal T$ then
$H$ is freely generated by the non-trivial elements
\begin{equation}
tx(\phi(tx))^{-1},
\label{eqS10}
\end{equation}
where $t$ ranges over $\mathcal T$ and $x$ ranges over $X$ (see \cite{Lyn}).
This set is called a Schreier basis for $H$.
Since $tx = \phi(tx)$ only
when $tx \in \mathcal T$ then by (\ref{eqS10}) the rank of $H$ equals
the {\em cyclomatic number} of $\mathcal G$,
the cardinality of the ``missing'' edges in the directions of $X$
in $\mathcal T$, that is
\begin{equation}
\mbox{rank}(H) = |\{ e \in E({\mathcal G}) - E({\mathcal T}) \}|,
\label{eqS20}
\end{equation}
where $E({\mathcal G})$, $E({\mathcal T})$ denote the set of edges of
$\mathcal G$, $\mathcal T$ respectively.
This is because each edge is labeled with some $x \in X$ in exactly
one direction, and thus counted exactly once.

Suppose now that $F$ is finitely generated with $\mbox{rank}(F)=n$,
and $H$ is of finite index $m$ in $F$. Then $|E({\mathcal G})|=nm$ and
$|E({\mathcal T})|=m-1$ (since $\mathcal T$ is a tree).
By (\ref{eqS20})
we get that
\begin{equation}
\mbox{rank}(H) = 1 + (n-1)m.
\label{eqS23}
\end{equation}
This is  Schreier Formula (see \cite{Lyn}). When $H$ is not
necessarily of finite index in $F$ and also not necessarily finitely
generated, we give in the proposition
below a formula that generalizes the above one. The rank is
computed on a Schreier transversal, and the simpler form of the
formula is given in Corollary~\ref{crS10}, which expresses the
rank in terms of the cogrowth (see below) of the subgroup.
The common way of computing the rank of the subgroup as a
limit of the ranks of the fundamental groups (the cyclomatic numbers)
of finite subgraphs deals with counting {\em edges}. Whereas, what we
are doing here is counting only {\em vertices}. \\

We use the following terminology and notation on graphs. A {\em path}
in a graph $\mathcal G$ is a sequence $v_0,e_1,v_1,e_2, \ldots$,
$v_i \in V({\mathcal G})$, $e_i \in E({\mathcal G})$, such that $e_i$
starts at the vertex $v_{i-1}$ and terminates at $v_i$. The
length of a path $v_0,e_1,v_1,e_2, \ldots, v_n$ is $n$.
A {\em simple path} is a path in which the vertices along it are distinct,
except possibly for the first and last one, in which case it is a
{\em simple closed path} or a {\em simple circuit}.
We assume that each path is {\em reduced}, i.e. it is not homotopic
to a shorter one when the initial and terminal vertices are kept fixed.

If ${\mathcal H} \subseteq {\mathcal G}$, i.e.
$\mathcal H$ is a collection of vertices
and edges of the graph $\mathcal G$, then
we denote by $< {\mathcal H} >$ the subgraph {\em generated} by $\mathcal H$.
It is the smallest
subgraph of $\mathcal G$ which contains $\mathcal H$.
That is, we add to $\mathcal H$ the
endpoint vertices of all the edges in $\mathcal H$. On the other hand, the
subgraph of $\mathcal G$ {\em induced} by $\mathcal H$ is the one whose
vertices are
those of $\mathcal H$ and whose edges are all the edges which join
these vertices
in $\mathcal G$. An induced subgraph is a subgraph which is induced by some
${\mathcal H} \subseteq {\mathcal G}$. If ${\mathcal H}_1, {\mathcal H}_2
\subseteq \mathcal G$
then ${\mathcal H}_1 - {\mathcal H}_2$ is the collection of
vertices $V({\mathcal H}_1) - V({\mathcal H}_2)$ and edges $E({\mathcal H}_1)
- E({\mathcal H}_2)$,
and it does not necessarily form a subgraph of $\mathcal G$,
even when ${\mathcal H}_1$
and ${\mathcal H}_2$ are subgraphs of $\mathcal G$. The {\em boundary}
of the subgraph $\mathcal H$ of $\mathcal G$ is
$\partial {\mathcal H} =  {\mathcal H} \cap <{\mathcal G} - {\mathcal H}>$,
and its {\em interior} is $\dot{\mathcal  H} = {\mathcal H} -
\partial {\mathcal H}$.
The {\em outer boundary} of $\mathcal H$ (in $\mathcal G$)
is the set of vertices
of ${\mathcal G} - {\mathcal H}$ which are adjacent to $\mathcal H$
in $\mathcal G$.
Assume now that each edge of $\mathcal G$ is labeled with some $x \in X$
in one direction and with $x^{-1} \in X^{-1}$ in the other direction.
Then we define $E^X_{out}({\mathcal H})$ to be the set of edges of
${\mathcal G} - {\mathcal H}$ whose initial vertices with respect to
the directions $X$ are in $\mathcal H$.
If ${\mathcal H}_i \subseteq {\mathcal G}$, $i=1,2, \ldots$, then $\mathcal H =
\liminf {\mathcal H}_i$ if the vertices of $\mathcal H$ are
$V({\mathcal H}) = \bigcup_{i \geq 1} \bigcap_{j \geq i} V({\mathcal H}_j)$,
and its edges are
$E({\mathcal H}) = \bigcup_{i \geq 1} \bigcap_{j \geq i} E({\mathcal H}_j)$.

Finally, let
$\alpha({\mathcal H})=
|\pi_0({\mathcal H})|$ be the cardinality of the (connected)
components of $\mathcal H$.
\begin{proposition}
Let $F$ be a free group of rank $n$ and let $H < F$. Let $\mathcal T$ be
a Schreier transversal for $H$ in $F$ and let ${\mathcal T}_i$ be finite
subgraphs of $\mathcal T$ such that ${\mathcal T} = \liminf {\mathcal T}_i$.
Then
\begin{equation}
\mbox{\em rank}(H) = \lim_{i \rightarrow \infty}
\left(
\alpha({\mathcal T}_i) + (n - 1)|V({\mathcal T}_i)| - \frac{1}{2}
\sum_{j=1}^{\alpha({\mathcal T}_i)} |V(\partial_{out}{\mathcal T}_{i,j})|
\right),
\label{eqS22}
\end{equation}
where, for a fixed $i$, $\partial_{out} {\mathcal T}_{i,j}$ is the outer
boundary (in $\mathcal T$) of the component ${\mathcal T}_{i,j}$
of ${\mathcal T}_i$, for $j= 1, \ldots, \alpha({\mathcal T}_i)$.
\label{prS40}
\end{proposition}
{\em Proof}.
If $H$ is of finite index $m$ in $F$ then there exists $i_0$
such that ${\mathcal T}_{i} = \mathcal T$ for every $i \geq i_0$, and
then $\alpha({\mathcal T}_i) = 1$, $|V({\mathcal T}_i)| = m$ and
$|V(\partial_{out}{\mathcal T}_i)| = 0$. Thus (\ref{eqS22}) reduces
to Schreier Formula.

Assume that $H$ is finitely generated but of infinite index.
Denote as before by $\mathcal G$ the cosets graph of $H$, which contains the
Schreier transversal tree $\mathcal T$. Let $C({\mathcal G})$ be the {\em core}
of $\mathcal G$ (see \cite{Sta}), that is the minimal deformation retract
of $\mathcal G$. It is the minimal connected subgraph of $\mathcal G$ which
contains all its simple circuits.
Since $H$ is finitely generated $C({\mathcal G})$ is
finite, and there exists $i_0$ such that, after possibly renaming the
components of each ${\mathcal T}_i$, $V(C({\mathcal G}))$ is contained in
$V({\mathcal T}_{i,1})$ for each $i \geq i_0$. Let us denote by
${\mathcal G}_{i,j}$ the subgraph of $\mathcal G$ induced by ${\mathcal T}_{i,j}$.
Then $E({\mathcal G}) - E({\mathcal T}) = E({\mathcal G}_{i,1}) -
E({\mathcal T}_{i,1})$, for each $i \geq i_0$, and by (\ref{eqS20}) the
cardinality of this set equals the rank of $H$. Hence it suffices to
show that for each $i \geq i_0$ and for each $j$
\begin{eqnarray}
|E({\mathcal G}_{i,j}) - E({\mathcal T}_{i,j})| &=&
|E({\mathcal G}_{i,j})| - |E({\mathcal T}_{i,j})|
\label{eqS27} \\
&=& 1+(n - 1)|V({\mathcal T}_{i,j})|
- \frac{1}{2} |V(\partial_{out}{\mathcal T}_{i,j})|.
\label{eqS28}
\end{eqnarray}
So assume $i \geq i_0$. Then $|E^X_{out}({\mathcal G}_{i,j})| =
|E^X_{out}({\mathcal T}_{i,j})| =
|V(\partial_{out}{\mathcal T}_{i,j})|/2$ for each $j$,
since all simple circuits of $\mathcal G$ are in ${\mathcal G}_{i,1}$.
Each vertex in $\mathcal G$ is the initial vertex of exactly $n$ edges
in the directions $X$. Therefore
\begin{equation}
|E({\mathcal G}_{i,j})|= n|V({\mathcal G}_{i,j})|-|E^X_{out}
({\mathcal G}_{i,j})|.
\end{equation}
As for ${\mathcal T}_{i,j}$, since it is a tree then
\begin{equation}
|E({\mathcal T}_{i,j})| = |V({\mathcal T}_{i,j})| -1.
\end{equation}
Substituting in (\ref{eqS27}) gives (\ref{eqS28}).

Assume now that $H$ is not finitely generated. Then, because in general
\begin{equation}
|E^X_{out}({\mathcal G}_{i,j})| \geq
|V(\partial_{out}{\mathcal T}_{i,j})|/2,
\end{equation}
we get that for each $i,j$
\begin{equation}
|E({\mathcal G}_{i,j}) - E({\mathcal T}_{i,j})| \leq
1+(n - 1) |V({\mathcal T}_{i,j})| - \frac{1}{2}
|V(\partial_{out}{\mathcal T}_{i,j})|.
\end{equation}
Since $\mbox{rank}(H) = \lim_{i \rightarrow \infty}
\left(
\sum_{j=1}^{\alpha({\mathcal T}_i)} |E({\mathcal G}_{i,j}) -
E({\mathcal T}_{i,j})|
\right) = \infty$,
equation~(\ref{eqS22}) follows.
\hfill $\Box$\\

We remark that instead of taking finite subgraphs ${\mathcal T}_i$ such
that ${\mathcal T} = \liminf {\mathcal T}_i$,
the rank formula can be clearly given
as the supremum, over all finite subgraphs of $\mathcal T$, of the expression
appearing in (\ref{eqS22}).

A special case of Proposition~\ref{prS40} is when each component
${\mathcal T}_{i,j}$ is  a {\em ball}, that is its vertices are all the
vertices of $\mathcal T$ which lie at distance not greater than some fixed
$k$ from some fixed vertex. If $\mathcal H$ is a subgraph of $\mathcal G$ and
$|V({\mathcal G})| >1$ then we define $\delta({\mathcal H})$ to be the
number of components of $\mathcal H$ which consist of a single vertex,
i.e. balls of radius $0$. When $|V({\mathcal G})| =1$ then
$\delta({\mathcal G})$
is defined to be $0$.

When the ${\mathcal T}_i$ are concentric balls centered at the identity
$1$ then the values $|V({\mathcal T}_i)|$,
$i=0,1,2,\ldots$ relate to the {\em growth} function $\Gamma_{\dt}$
of $\mathcal T$, as is defined below.
By $l(g)$ we denote the {\em length} of $g \in F$, and we always assume
that the group elements are written in reduced form with respect to
the generating set $X$ of $F$. Then define
\begin{eqnarray}
\gamma_{\dt}(i) &=& |\{ v \in {\mathcal T} \mid l(v)=i \}|, \\
\Gamma_{\dt}(i) &=& |\{ v \in {\mathcal T} \mid l(v) \leq i \}|.
\end{eqnarray}
When $\mathcal T$ is a {\em minimal} Schreier transversal tree, that is
when it has also the property that every coset of $H$ is
represented by an element of minimal length, then $\Gamma_{\dt}(i)$
is the {\em cogrowth} function of $H$, relative to the generating
set of $F$, and is denoted by $\Gamma_{\dfh}(i)$ (see \cite{Ros1}).
We may look at $\Gamma_{\dfh}(i)$ as representing the ``volume''
of the ball of radius $i$ with center $1$ in the cosets graph of $H$
(with the metric induced by the word metric on $F$).
If, in addition, $H$ is a normal
subgroup of $F$ then the cogrowth function of $H$ equals the
growth function of the group $F/H$, relative to the the generating
set which is the canonical image of the generating set of $F$.
(In this case the Schreier transversal for $H$ which is minimal with
regard to a fixed ShortLex order on $F$ is also suffix-closed.)
\begin{corollary}
Let $F$ be a free group of rank $n$, let $H$ be a subgroup of $F$ and
let $\mathcal T$ be
a Schreier transversal for $H$ in $F$. Let ${\mathcal T}_i$ be induced
finite subgraphs of $\mathcal T$, whose components ${\mathcal T}_{i,j}$ are
balls, such that ${\mathcal T} = \liminf {\mathcal T}_i$. Then
\begin{eqnarray}
\mbox{\em rank}(H) &=& \lim_{i \rightarrow \infty}
\left(
\alpha({\mathcal T}_i)  - \frac{1}{2} \delta({\mathcal T}_i)
+ (n - 1)|V({\mathcal T}_i)|
- \frac{2n-1}{2} |V(\partial {\mathcal T}_i)|
\right) \\
&=& \lim_{i \rightarrow \infty}
\left(
\alpha({\mathcal T}_i) - \frac{1}{2} \delta({\mathcal T}_i)
+ (n - 1)|V(\dot{\mathcal T}_i)|
- \frac{1}{2} |V(\partial {\mathcal T}_i)|
\right).
\end{eqnarray}
In particular,
\begin{eqnarray}
\mbox{\em rank}(H) &=& 1 + \lim_{i \rightarrow \infty}
\left(
(n-1)\Gamma_{\dt}(i) - \frac{1}{2} \gamma_{\dt}(i+1)
\right) \\
&=& 1 + \lim_{i \rightarrow \infty}
\left(
(n-1)\Gamma_{\dfh}(i) - \frac{1}{2} \gamma_{\dfh}(i+1)
\right).
\label{eqS35}
\end{eqnarray}
\label{crS10}
\end{corollary}
{\em Proof}.
The corollary follows from the fact that when the core $C({\mathcal G})$ is
finite then for each $i$ large enough every vertex
of $\partial {\mathcal T}_{i,j}$ is adjacent to $2n-1$ vertices of ${\mathcal T}
- {\mathcal T}_{i,j}$, unless ${\mathcal T}_{i,j}$
is a single vertex and then it
is adjacent to $2n$ vertices of ${\mathcal T} - {\mathcal T}_{i,j}$.
When $H$ is not finitely generated then we first notice that the
expression we calculate for each ball is non-negative. Secondly,
since ${\mathcal T} = \liminf {\mathcal T}_i$, then for every $r$ there exists
$i_r$ such that, for every $i \geq i_r$, ${\mathcal T}_i$ has a component
(ball) which contains the ball of radius $r$ around the identity.
But the expression calculated on these balls tends to infinity whenever
$H$ is of infinite rank, as shown below. This can also be concluded
directly from Proposition~\ref{prS40}.
\hfill $\Box$
\section{Rank-growth}
Given a Schreier transversal $\mathcal T$, let us define
\begin{eqnarray}
r_{\dt}(i) &=&
1+(n-1)\Gamma_{\dt}(i) - \frac{1}{2} \gamma_{\dt}(i+1) \\
&=& 1+\frac{2n-1}{2}\Gamma_{\dt}(i) - \frac{1}{2}\Gamma_{\dt}(i+1).
\label{eqS42}
\end{eqnarray}
$r_{\dt}(i)$ is an upper bound to the cyclomatic number of the subgraph
of $\mathcal G$ which is induced by the vertices of $\mathcal T$ of distance
at most $i$ from the root.
In case $\mathcal T$ is a minimal Schreier transversal then $r_{\dt}(i)$
is also denoted by $r_{\dh}(i)$:
\begin{equation}
r_{\dh}(i) = 1+(n-1)\Gamma_{\dfh}(i)-\frac{1}{2} \gamma_{\dfh}(i+1).
\label{eqS44}
\end{equation}
The sequence $r_{\dt}(i)$, $i=1,2, \ldots$ is non-decreasing. This is because
\begin{eqnarray}
r_{\dt}(i) - r_{\dt}(i-1) &=&
\frac{2n-1}{2}\gamma_{\dt}(i) - \frac{1}{2}\gamma_{\dt}(i+1),
\end{eqnarray}
and each vertex of $\mathcal T$ of level $i$ is adjacent to at most $2n-1$
vertices of level $i+1$. Thus $r_{\dt}(i)$ becomes eventually
constant if and only if either $\mathcal T$ is finite, or for some $i_0$
each vertex of $\mathcal T$ of level $i \geq i_0$ has
degree exactly $2n$, and this happens if and only if there are only
finitely many edges in $E({\mathcal G}) - E({\mathcal T})$, or equivalently when
$H$ is finitely generated.

It is interesting to know also the
rate in which the function $r_{\dt}(i)$ grows. A preorder is defined
on growth functions by
\begin{equation}
f_1(i) \preceq f_2(i) \ \Llrw \
\exists c>0 \ \forall i \ [ f_1(i) \leq c f_2(ci) \ ].
\end{equation}
Then an equivalence relation is given by
\begin{equation}
f_1(i) \sim f_2(i) \ \Llrw \ f_1(i) \preceq f_2(i)
\ \mbox{and} \ f_2(i) \preceq f_1(i).
\end{equation}
(we refer to \cite{Gri} for a survey on growth functions of groups
and to Gromov's \cite{Gro} rich and beautiful geometric theory.)
In Theorem~\ref{thS60} below we show that when the subgroup $H$ of the
free group $F$ is {\em supnormal}, i.e. contains a non-trivial subgroup
which is normal in $F$, then
for every Schreier transversal $\mathcal T$
of $H$, its growth function $\Gamma_{\dt}(i)$ is equivalent to
the function $r_{\dt}(i)$.
This implies that the cogrowth of $H$
is also equivalent to what we call the rank-growth of $H$.
We look at $H$ as the direct limit of the subgroups
\begin{equation}
H_i = < \{h \in H \mid l(h) \leq i\} >,
\end{equation}
where $l(h)$ is measured with respect to the generating set of $F$.
Then the {\em rank-growth} of $H$ (with respect to the generators of $F$) is
\begin{equation}
rk_{\dh} (i) = \mbox{rank} (H_i).
\end{equation}
Clearly, if we choose another generating set for $F$, we get an equivalent
rank-growth function. Notice that $H_i$
is the fundamental group of the subgraph of the cosets graph $\mathcal G$ of
of $H$ which contains all paths starting at $1$ of length $\leq i$.
Thus $rk_{\dh} (i)$ is a non-decreasing function. If we define
\begin{equation}
\rho_{\dh} (i) = \mbox{rank} (\pi_1({\mathcal B}_i)),
\end{equation}
where ${\mathcal B}_i$ is (the induced subgraph which is) the ball of radius
$i$ centered at the vertex $1$
of $\mathcal G$, then
\begin{equation}
\rho_{\dh} (i) = rk_{\dh} (2i+1).
\end{equation}
Therefore $rk_{\dh}(i)$ and $\rho_{\dh}(i)$ are equivalent.
Also $\rho_{\dh}(i) \sim r_{\dh}(i)$. In fact,
\begin{equation}
\rho_{\dh}(i) \leq r_{\dh}(i) \leq \rho_{\dh}(i+1).
\end{equation}
More precisely,
\begin{equation}
r_{\dh}(i) = \rho_{\dh}(i) +\frac{1}{2}(|E^X_{out}({\mathcal B}_i)| -
\gamma_{\dfh}(i+1)) \leq \rho_{\dh}(i+1).
\label{eqS54}
\end{equation}
\begin{theorem}
Let $H$ be a supnormal subgroup of a finitely generated free group $F$,
and let $\mathcal T$ be a Schreier transversal for $H$ in $F$. Then
\begin{equation}
r_{\dt}(i) \sim \Gamma_{\dt}(i).
\end{equation}
In fact, if $H$ is not necessarily supnormal but has the property that
$|F:N_{\df}(A)| < \infty$ for some non-trivial $A<H$ then
\begin{equation}
rk_{\dh}(i) \sim \Gamma_{\dfh}(i).
\end{equation}
\label{thS60}
\end{theorem}
{\em Proof}.
For every Schreier transversal of a subgroup of $F$ we have
$r_{\dt}(i) \preceq \Gamma_{\dt}(i)$. This follows immediately from
the definition of $r_{\dt}(i)$ - see (\ref{eqS42}).

Suppose now that $H$ is supnormal. Let $h$ be a non-trivial
element of a subgroup of $H$ which is normal in $F$, and let $m = l(h)$
(as usual, the length is with respect to the generators of $F$).
Then at every vertex $v$ of the cosets graph $\mathcal G$ of $H$, if we follow
the path defined by $h$ we form a circuit. Therefore at every vertex
of $\mathcal T$ of level at most $i$, by following the path defined by $h$
we reach a vertex of $\mathcal T$ of level at most $i+m$ where we must
stop because the next edge is missing. The number of these missing edges
is less then or equal to $r_{\dt}(i+m)$. Since at most
$m$ vertices are the starting point of a tour defined by $h$ which reaches
the same missing edge $h$ then
\begin{equation}
\Gamma_{\dt}(i) \leq m r_{\dt}(i+m).
\end{equation}
By the two inequalities we have
\begin{equation}
r_{\dt}(i) \sim \Gamma_{\dt}(i).
\end{equation}
Applying this result to a minimal Schreier transversal yields
\begin{equation}
rk_{\dh}(i) \sim r_{\dh}(i) \sim \Gamma_{\dfh}(i).
\end{equation}

The condition of $H$ being supnormal can be weakened. It suffices to
demand
that $H$ contains a non-trivial subgroup $A$ such that $|F:N_{\df}(A)|
< \infty$,
because then the cogrowth of $H$ is equivalent to the growth (with respect
to the generators of $F$) of the minimal coset representatives of
$H \cap N_{\df}(A)$ in $N_{\df}(A)$. Even more, we need only
the growth (again, with respect to the generators of $F$) of
$\{ g \in {\mathcal T} \mid gAg^{-1} \subseteq H \}$, where $\mathcal T$
is a minimal Schreier transversal for $H$,
to be equivalent to the cogrowth of $H$ in $F$.
\hfill $\Box$ \\

Since $\Gamma_{\dt}(i) \preceq \Gamma_{\dfh}(i)$ for every Schreier
transversal $\mathcal T$ of a subgroup $H$ of $F$, then by Theorem~\ref{thS60}
when $H$ is supnormal in $F$ then $r_{\dt}(i) \preceq rk_{\dh}(i)$.
We also notice that a special case of Theorem~\ref{thS60} is the known result
stating that a supnormal subgroup of a finitely generated group is of
finite index if and only if it is finitely generated. And when $H$ is
normal in $F$, then the growth $\Gamma_{\dg}(i)$ of the group $G=F/H$
is equivalent to the rank-growth of $H$ and to the growth of
\begin{equation}
r_{\dh}(i) =1+(n-1)\Gamma_{\dg}(i) - \frac{1}{2} \gamma_{\dg}(i+1).
\end{equation}

The growth of the subgroup $H$ is always exponential when it is of
rank greater than $1$, since it is free. But Grigorchuk showed
(\cite{Gri0}) that when $H$ is normal then its ``growth exponent''
$\limsup_{i \rightarrow \infty} \Gamma_{\dh}^{(\df)} (i)^{1/i} = 2n-1$,
if and only if $G=F/H$ is amenable,
(in fact, Grigorchuk \cite{Gri0} obtained more: a formula which connects the
growth exponent of $G$ with the spectral radius of a random walk on $G$),
where $n = \mbox{rank}(F)$ and
$\Gamma_{\dh}^{(\df)} (i)$ represents the growth of $H$ with respect
to the generators of $F$. (Recall that a group $G$ is amenable if there exists
an invariant mean on $B(G)$, the space of all bounded complex-valued
functions on $G$ with the sup norm $\parallel f \parallel_{\infty}$, see
\cite{Gre}). When $G$ is non-amenable then the growth
exponent of $H$ is less than $2n-1$. But then the group $G$ has exponential
growth, and we have shown that in this case the rank-growth of $H$ is
also exponential, i.e. the maximal possible (up to equivalence).
This seems at first sight contradictory. To illustrate this phenomenon
we may think of a tree, called $F$, that we prune its sides going from
bottom upward. The number of branches we cut is called (half) the rank
of $H$, the tree that is left after the pruning is called $G$, and
(part of) what we cut is called $H$. Then
the further the cut is from the periphery and closer to the middle
of the tree the larger $H$ is, the smaller $G$ is, and the rank of $H$
also becomes smaller since we cut towards the main branches.

Although the rank of the subgroup of a free group can be expressed, as
we have seen in Corollary~\ref{crS10}, in terms of the growth function of
any Schreier transversal of it, the growth function itself of one
Schreier transversal of an infinitely generated subgroup may in general
differ completely from that of another Schreier transversal.
This is shown in the next proposition.
\begin{proposition}
There exists a subgroup of the free group of rank $2$ with
exponential cogrowth which has a Schreier transversal $\mathcal T$
whose growth is $\Gamma_{\dt}(i) = i+1$.
\label{prS50}
\end{proposition}
{\em Proof}.
We will construct the cosets graph of such a subgroup inductively.
Let $c$ be a positive integer which is large enough. First we make a
simple circuit $\lambda_1$ of length $c$ that starts at the root $1$.
Then at the $n$-th step we construct a path $\lambda_n$ of length $2nc$,
whose vertices, apart from the initial and terminal ones, are new.
The initial vertex of $\lambda_n$ is the one before the last vertex
in the path $\lambda_{n-1}$. The terminal vertex of $\lambda_n$ is
chosen to be of minimal distance from the root among the vertices
whose degree is less than $4$.

If we delete the last edge of each path $\lambda_n$, then we get a
linear Schreier transversal $\mathcal T$, i.e.  $\Gamma_{\dt}(n) = n+1$.
On the other hand, if we delete the middle edge of each $\lambda_n$,
then we get a tree ${\mathcal T}'$ with exponential growth,
because each vertex of it has degree $4$, except for a sequence of
vertices $v_n$ of distances $\geq cn$ respectively from the root. Since
the cogrowth function is greater than or equal to the growth function of any
Schreier transversal of the subgroup, the result follows.
\hfill $\Box$ \\

It is shown in \cite{Ros1} that when $H = H_1 \cap H_2$ the
cogrowth of $H$ satisfies
\begin{equation}
\Gamma_{\dfh}(i) \leq \Gamma_{\dfhA}(i) \Gamma_{\dfhB}(i) \ \ \
\mbox{for every $i$}.
\end{equation}
The rank-growth of the intersection of two subgroups behaves similarly.
\begin{proposition}
Let $H_1, H_2$ be non-trivial subgroups of a finitely generated
free group $F$ and let $H = H_1 \cap H_2$. Then
\begin{equation}
rk_{\dh}(i) \leq 1 + 2(rk_{\dhA}(i)-1)(rk_{\dhB}(i)-1)-\min(rk_{\dhA}(i),
rk_{\dhB}(i))
\end{equation}
for every $i$.
Hence
\begin{equation}
rk_{\dh}(i) \preceq rk_{\dhA}(i)rk_{\dhB}(i).
\end{equation}
\label{prS60}
\end{proposition}
{\em Proof}.
This follows immediately from the best general bound for the intersection
of finitely generated subgroups in free groups, which is due to
Burns (\cite{Bur}).
\hfill $\Box$
\section{The Generalized Word Problem}
Given a subgroup $H$ of a group $G$ it is interesting to know the
{\em distorsion} of $H$ with respect to $G$, that is a bound $f(i)$ of
the length, with respect to a finite set of generators of $H$, of an
element of $H$ whose length is $i$ with respect to a finite set of
generators of $G$ (see \cite{Gro}, \cite{Farb}).
When $F$ is free then it is known that every element of a subgroup $H$
of it, whose length is $i$ in $F$, has length at most $i$ with respect to a
Schreier basis of $H$ (or a Nielsen-reduced basis, which is no other than a
minimal Schreier basis), thus the distortion
is linear. A bit more precise description is obtained by using
$d(w,{\mathcal T})$, the distance of $w \in F$ from a Schreier
transversal $\mathcal T$. That is
\begin{equation}
d(w,{\mathcal T}) = \min \{ l(t^{-1}w) \mid t \in {\mathcal T} \},
\end{equation}
where $l$ denotes the length in $F$. Notice that $d(w,{\mathcal T}) \leq l(w)$
since $1 \in {\mathcal T}$. Then if $B_{\dt}$ is a corresponding
Schreier basis for the subgroup $H<F$ then every
$w \in F$ can be written in the form
\begin{equation}
w = b_{i_1}^{\varepsilon_1} \cdots b_{i_k}^{\varepsilon_k} \bar{w},
\end{equation}
with $b_{i_j} \in B_{\dt}$, $\varepsilon_j= \pm 1$ and $\bar{w} \in
\mathcal T$, such that $k \leq d(w,{\mathcal T})$.
To see it,
let $t \in {\mathcal T}$ be the maximal prefix of $w$ in $\mathcal T$, i.e.
$l(t^{-1}w) = d(w,{\mathcal T})$. If $t=w$ we are done. Otherwise, there
exists $x \in X \cup X^{-1}$, $X$ the generating set of $F$, such that
$b_{i_1}^{\varepsilon_1}=tx(\phi(tx))^{-1} \in B_{\dt} \cup B_{\dt}^{-1}$,
$\phi$ the coset map, and such that $w=txw'$ when written in reduced
form. Thus
\begin{equation}
w = b_{i_1}^{\varepsilon_1}\phi(tx)w'.
\end{equation}
But $d(\phi(tx)w',{\mathcal T}) \leq l(w') = d(w,{\mathcal T})-1$
and we proceed by induction.

The above shows that when we are given a
Schreier transversal $\mathcal T$ and a corresponding Schreier basis
$B_{\dt}$ for $H < F$ then it is possible to obtain algorithmically
a ``normal form'' modulo $H$ for every element of $F$, i.e. its coset
representative in $\mathcal T$, and this demonstrates the importance of
Schreier generating sets (whose shape and role is similar to those of
Gr\"obner bases for algebras, see \cite{Ros2}). Thus the generalized
word problem
for $H$ in $F$ is then solvable. In fact, whenever $G = <X \mid R>$,
and $H$ is a subgroup of $G$ generated by a set $S \subseteq F=<X>$,
then the generalized word problem for $H$ in $G$ is solvable when
the underlying set of the subgroup $<N,S>$ of $F$, where $N=<R>^{\df}$
is the normal closure of $R$ in $F$, as well as a set $\mathcal T$ of coset
representatives for $H$ in $G$, are recursively enumerable (r.e.) sets.
For, by listing the elements of $<N,S>$ and those of $\mathcal T$ we can list
all elements of $F$, and also find which coset is the coset $<N,S>$ in $F$.
Hence, both the set $<N,S>$ and its complement in $F$ are r.e. and
therefore $<N,S>$ is recursive.
\begin{proposition}
Let $G = <X \mid R>$ be a finitely generated group and let $H$ be a
subgroup of $G$ generated by a set $S$. Suppose also that the subgroup
$A=<N,S>$ of the free group $F=<X>$, where $N=<R>^{\df}$, is r.e.
Then the generalized word problem for $H$ in $G$ is solvable
whenever one of the functions $\Gamma_{\dfa}(i), r_{\da}(i)$ or
$rk_{\da}(i)$ is recursive.
\label{prS80}
\end{proposition}
{\em Proof}.
First we remark that $A$ is r.e. for example when $R$ and $S$ are r.e.
We construct inductively ${\mathcal B}_i$, the concentric
balls of radius $i$, of the cosets graph of $A$ in $F$.
We start with the vertex $1$. Then assuming that ${\mathcal B}_i$ was constructed,
we first extend it to level $i+1$ without forming new circuits. If the
number of vertices at level $i+1$ agrees with $\Gamma_{\dfa}(i+1)$ or
with $r_{\da}(i+1)$ or if $\pi_1({\mathcal B}_{i+1}) = rk_{\da}(2i+1)$ then
we are done. Otherwise, by listing the elements of $A$, each defining
a circuit in the cosets graph, we stop until we reach the desired values
of our functions.
\hfill $\Box$ \\ \\
\section{Contractable Spanning Subcomplexes}
If we look at Proposition~\ref{prS40} we see that it makes little use
of the group structure. It is mainly a statement about $m$-regular
graphs, i.e. graphs whose vertices have all the same degree $m$.
We may then try to generalize this theorem from such graphs to special
simplicial complexes. When $\CC$ is a simplicial complex
then we denote by $| {\CC} |$ the topological space of $\CC$ (as in
\cite{Spa}), and whenever we relate some topological properties to
$\CC$ they describe, in fact, those of $| {\CC} |$.
We Call a subcomplex $\DD$ of a $d$-dimensional simplicial complex
$\CC$ a {\em spanning} subcomplex if
\begin{description}
\item[(i)] $\DD$ contains ${\CC}^{(d-1)}$, the $(d-1)$-skeleton of $\CC$;
\item[(ii)] every principal simplex of $\DD$ (i.e. a simplex which is
not a face of another simplex of $\DD$ of higher dimension) is also
principal in $\CC$.
\end{description}

Then the analogue of a spanning tree in graph theory
is a spanning subcomplex $\DD$ whose topological space is contractible.
We call such a subcomplex a {\em contractible spanning subcomplex}.
In this case condition (ii) becomes redundant. This is because if
$\sigma_1^{d-1}$ is principal in $\DD$ but is a face of a $d$-simplex
$\sigma_2^d$ of $\CC$ then by considering the boundary of $| \sigma_2^d |$
in $| {\DD} |$ we get that $\pi_{d-1} (| {\DD} |)$ is not trivial and thus
$\DD$ is not contractible. Or we can look at the homology of $\DD$,
and see that $H_{d-1}({\DD})$ does not vanish since $\sigma_1^{d-1}$
does not appear in $B_{d-1}({\DD})$ but is a summand of a cycle.

We come now to the analogue of Proposition~\ref{prS40} for simplicial
complexes. The formula we give, however, is not so nice as the one
in the one-dimensional case, where all terms involve only the
zero-dimensional skeleton. We use the following additional
notation and definitions.
Let $\CC$ be a simplicial complex and let $\DD$ be a subcomplex of it.
The collection of $k$-simplices of $\CC$ is denoted by $F^k({\CC})$,
and its cardinality is denoted by $\beta_k({\CC})$. When $X$ is a collection of
simplices of $\CC$ we denote by $< X >$ the subcomplex generated by
$X$. If ${\DD}_1, {\DD}_2$ are subcomplexes of $\CC$ then ${\DD}_1 - {\DD}_2$
is the collection of simplices ${\DD}_1 - {\DD}_2 = \{ \sigma \mid \sigma
\in {\DD}_1, \ \sigma \notin {\DD}_2 \}$, and it does not necessarily form a
subcomplex. The {\em boundary} of the subcomplex $\DD$ of $\CC$ is
$\partial \DD =  \DD \cap <{\CC} - {\DD}>$, and its {\em interior} is
$\dot{\DD} = {\DD} - \partial {\DD}$.
If $\sigma^k$ is a $k$-simplex of $\CC$ we define its {\em degree}
to be $\mbox{deg}(\sigma^k)=|\{\sigma^{k+1} \in
F^{k+1}({\CC}) \mid \sigma^k \subset \sigma^{k+1} \}|$. If $X$ is a
collection of simplices of $\CC$ then $\mbox{deg}(X)$ is the total
degree of the members of $X$. When all members of $X$ belong to a
subcomplex $\DD$ of $\CC$ and we want to compute the degree relative
to this subcomplex then we write it $\mbox{deg}_{\dd}(X)$.
When $\CC$ is of dimension $d$ then we say it is
$n$-{\em regular} if every $(d-1)$-simplex of it has degree $n$.
If ${\DD}_i$, $i \geq 1$, is a sequence of subcomplexes of $\CC$ then
we denote by $\liminf {\DD}_i$ the subcomplex of $\CC$ whose simplices
are $F^k(\liminf {\DD}_i) = \bigcup_{i \geq 1}
\bigcap_{j \geq i} F^k({\DD}_j)$, for every $k \geq 0$.
\begin{theorem}
Let $\CC$ be a countable $d$-dimensional simplicial complex
which contains a contractible spanning subcomplex $\DD$. Then
$| {\CC}  |$ is homotopic to a bouquet of $r$ ($r$ can be $\infty$)
$d$-spheres. \\
If, in addition, $\CC$ is $n$-regular and
${\DD}_i$ are finite subcomplexes of $\DD$ such that $\DD = \liminf {\DD}_i$
then
\begin{equation}
r = \lim_{i \rightarrow \infty}
\left(
\frac{1}{d+1} (n \beta_{d-1}({\DD}_{i}) - \mbox{\em deg}_{<\dd-\dd_i>}(F^{d-1}
(\partial {\DD}_{i}))) - \beta_d({\DD}_{i}) \right).
\label{eqCSS15}
\end{equation}
\label{thCSS20}
\end{theorem}
{\em Proof}.
We define a contractible space to be homotopic to a bouquet of zero
$d$-spheres. So assume that $| {\CC} |$ is not contractible.
Since $| {\DD} |$ is contractible, $| {\CC} |$ is homotopic to $| {\CC} | /
| {\DD} |$.
If there are $r$ ($r$ can be $\infty$) $d$-simplices $\sigma^d$ which
are not in $\DD$, then since the boundary of $| \sigma^d |$ is in
$| {\DD} |$, we get that $| {\CC} |$ is homotopic to a bouquet of $r$
$d$-spheres (see also \cite{Bjo} for shellable complexes, where every
shellable complex contains a contractible spanning subcomplex but not
necessarily the other way round).

The second part of the proof is similar to that of Proposition~\ref{prS40}.
In fact we are only
dealing with the simplices of $\CC$ of dimensions $d-1$ and $d$,
and then we compute the number of $d$-simplices of ${\CC} - {\DD}$.
\hfill $\Box$ \\

We remark that in case each subcomplex ${\DD}_i$ in (\ref{eqCSS15}) has
contractible
connected components then $\beta_d({\DD}_i)$ may be expressed
in terms of the $\beta_j({\DD}_i)$, $j=0, \ldots, d-1$ using
the (topological) {\em Euler characteristic}
\begin{equation}
\chi({\DD}_i^{(d-1)}) = \sum_{i=0}^{d-1} (-1)^i \beta_i({\DD}_i^{(d-1)}).
\end{equation}

\end{document}